\documentclass[11pt, a4paper]{amsart}

\usepackage{amscd}
\usepackage{a4}
\usepackage{amsmath}
\usepackage{amsxtra}
\usepackage[latin1]{inputenc}
\usepackage{graphicx}
\usepackage{amssymb}
\usepackage{pictex}
\usepackage{color}

\newtheorem{thm}{Theorem}[section]
\newtheorem{lemma}[thm]{Lemma}
\newtheorem{prop}[thm]{Proposition}
\newtheorem{cor}[thm]{Corollary}
\newtheorem{conj}[thm]{Conjecture}

\theoremstyle{definition}

\newtheorem{dfn}[thm]{Definition}

\newtheorem{rmk}[thm]{Remark}
\newtheorem{example}[thm]{Example}

\newcommand{\up}{\overline{d}}
\newcommand{\down}{\underline{d}}
\newcommand{\sspec}{_{\down, \up}}

\numberwithin{equation}{section}

\newcommand{\vecsp}{V\sspec}
\newcommand{\dimvecsp}{q}

\newcommand{\bettiset}{B\sspec}
\newcommand{\cbettiset}{\overline{B}\sspec}
\newcommand{\poset}{\varPi \sspec}
\newcommand{\pure}{\pi}
\newcommand{\cbetti}{\overline{\beta}}
\newcommand{\pdim}{p}
\newcommand{\MaxLex}{\beta^{\mathrm{lex}}}

\begin{document}

\title[Graded Betti Numbers of Cohen-Macaulay Modules] {Graded Betti
  Numbers of Cohen-Macaulay Modules and the Multiplicity Conjecture}

\author{Mats Boij} \address{Department of Mathematics, KTH \\ S--100
  44 Stockholm \\ Sweden} \email{boij@math.kth.se}

\author{Jonas Söderberg} \address{Department of Mathematics, KTH \\
  S--100 44 Stockholm \\ Sweden} \email{jonasso@math.kth.se}


\begin{abstract}
  We give conjectures on the possible graded Betti numbers of
  Cohen-Macaulay modules up to multiplication by positive rational
  numbers. The idea is that the Betti diagrams should be non-negative
  linear combinations of \emph{pure} diagrams. The conjectures are
  verified in the cases where the structure of resolutions are known,
  i.e., for modules of codimension two, for Gorenstein algebras of
  codimension three and for complete intersections.  The motivation
  for the conjectures comes from the Multiplicity conjecture of
  Herzog, Huneke and Srinivasan.
\end{abstract}

\maketitle


\section{introduction}
The graded Betti numbers of Cohen-Macaulay algebras and modules are
invariants that have been used in many different ways. In some special
cases, there are structure theorems that tells us precisely what the
graded minimal free resolutions can be, as in the case of the
Hilbert-Burch theorem for Cohen-Macaulay algebras of codimension two
or the Buchsbaum-Eisenbud theorem for Gorenstein algebras of
codimension three~\cite{Buchsbaum-Eisenbud}. There has been some
attempts to develop such structure theorems more generally
\cite{Kustin-Miller1, Kustin-Miller3}, but so far, there are
no more strong results in this direction. For some combinatorially
defined algebras, like the Stanley-Reisner rings of simplicial
polytopes, there are direct connections between the graded Betti
numbers and topological invariants.  In recent years, the Multiplicity
conjecture of Herzog, Huneke and Srinivasan~\cite{Herzog-Srinivasan, Huneke-Miller} 
has been inspiring a lot of interesting results
on the possible graded Betti numbers. The conjecture gives lower and
upper bounds for the multiplicity of a graded algebra in terms of the
lowest degrees and highest degrees where the Betti numbers are
non-zero. The conjecture has been verified in a number of interesting
cases, including the ones mentioned above where the structure of the
minimal free resolution is known (see Herzog and Srinivasan~\cite{Herzog-Srinivasan},
 Migliore, Nagel and R\"{o}mer~\cite{Migliore-Nagel-Roemer1, Migliore-Nagel-Roemer2}, 
Mir{\'o}-Roig~\cite{Miro-Roig} and
Kubitzke and Welker~\cite{Kubitzke-Welker}). However, the meaning of
the conjecture itself has been a bit mysterious, since the
multiplicity of an algebra is a very coarse invariant, while the
graded Betti numbers are much finer.
 
We will now argue for a stronger conjecture on the possible graded
Betti numbers of Cohen-Macaulay algebras and modules, which is easily
seen to imply the multiplicity conjecture in this case.  The
Multiplicity conjecture was inspired by the calculation of the
multiplicity, $e(A)$, of an algebra with a {\em pure resolution},
i.e., a resolution where each component is homogeneous, of degree
$d_i$, for $i=1,2,\dots,p$. This calculation shows that the
multiplicity is given by the very simple formula
$$
e(A) = \frac{d_1d_2\cdots d_p}{p!}.
$$
The Multiplicity conjecture gives the same expression as the lower and
upper bound for any graded algebra, where the degrees are replaced by
the lowest and highest degrees that appear in each component of the
minimal free resolution.
 
In our conjecture, Conjecture~\ref{MainConj}, the pure resolutions
also play a central role. We look at the possible Betti numbers of
graded Cohen-Macaulay algebras, or modules, as a subset of an integral
lattice, which sits inside a vector space over the rational numbers.
By taking direct sums of modules, we see that this subset is closed
under addition. If we also normalize so that the first Betti number,
$\beta_0$, is equal to one, the subset obtained will be closed under
convex combinations, i.e., linear combinations with non-negative
rational coefficients. The Betti numbers of the pure resolutions will
be in some sense extremal, and the idea is that they will actually be
the only vertices of the convex set described above. The pure
resolutions will be partially ordered in a natural way which is
compatible with the multiplicity given by the formula above. This
shows that the conjecture implies the multiplicity conjecture, as will
be made more precise in Proposition~\ref{ConjImplMC}.

In addition to this we conjecture (Conj.~\ref{ConvexConj}) that the
convex set also has a very nice structure of a simplicial polytope
which implies that the non-negative linear expansion of a Betti
diagram into pure diagrams is unique if the pure diagrams included
have to be totally ordered.  We are able to prove our conjectures in a
few, but interesting cases, where the set of possible Betti diagrams
is actually known. This includes the codimension two Cohen-Macaulay
modules generated in a single degree (Thm.~\ref{CodimOneTwo}),
codimension three Gorenstein algebras (Thm.~\ref{GorensteinCodim3Thm})
and complete intersections in any codimension (Thm.~\ref{CIThm}).

There is still much to be done in order to prove the conjectures, or
even to give more evidence for them, in a greater generality. For
example, in codimension three and higher, it is not known whether
there are modules with a pure resolution for each set of admissible
degrees, $d_1,d_2,\dots,d_p$. Furthermore, it would be much more
satisfactory to have proofs that use some kind of deformation
argument, showing that the modules actually split into sums of modules
with pure resolutions somewhere in the parameter space of
Cohen-Macaulay modules with fixed graded Betti numbers.


\section{The set of normalized Betti numbers}
Let $R = k[x_1, x_2, \dots, x_n]$ be the polynomial ring in $n$
variables over a field $k$. Consider $R$ as a graded ring by giving
each $x_i$ degree one.  All $R$-modules in this article are assumed to
be finitely generated and graded.  The $d$-th twist of an $R$-module
$M$, denoted by $M(d)$, is defined by $M(d)_i = M_{i+d}$.  If the
$R$-module $M$ has a minimal free resolution given by
$$
0 \to \bigoplus_j R(-j)^{\beta_{\pdim,j}} \to \dots \to \bigoplus_j
R(-j)^{\beta_{0,j}} \to M \to 0
$$
then $\beta_{i,j}(M) = \beta_{i,j}$ are the graded Betti numbers of
$M$ and the matrix $\beta(M)$ with entries $\beta(M)_{i,j} =
\beta_{i,j}(M)$ is called the \emph{Betti diagram} of $M$.  
From the
resolution we also see that the projective dimension of $M$ is $\pdim$
and when $M$ is Cohen-Macaulay it is equal to the codimension, that
is, $\dim R - \dim M = \pdim$.
Furthermore, if $M$ is Cohen-Macaulay
then the $h$-vector of $M$ is the Laurent polynomial $h_M(t) \in \mathbb{Z}[t,t^{-1}]$
given by
$$
h_M(t) = \frac{S_M(t)}{(1-t)^\pdim},
$$
where $S_M(t) = \sum_{i,j} (-1)^i \beta_{i,j}(M) t^j$. The multiplicity of 
$M$ can be calculated from its $h$-vector by
$$
e(M) = h_M(1).
$$
The fact that the polynomial
$S_M(t)$ is divisible by $(1-t)^\pdim$ is equivalent to the fact that
$$
\frac{d^i}{dt^i} S_M(1) = 0
$$
for $i = 0,1,\dots,\pdim-1$. This translates to the system of linear
equations
\begin{equation} \label{vSpaceEq} \left\{
    \begin{array}{lcc}
      \sum_{i,j} (-1)^i \beta_{i,j} &=& 0 \\
      \sum_{i,j} (-1)^i j\beta_{i,j} &=& 0 \\
      & \vdots & \\
      \sum_{i,j} (-1)^i j^{\pdim-1}\beta_{i,j} &=& 0
    \end{array}\right.
\end{equation}
whose set of solutions contains all possible Betti numbers of graded
Cohen-Macaulay modules of codimension $\pdim$.  This observation leads
to the following definition.

\begin{dfn}
  For any strictly increasing sequences of integers
  $\down=(\down_0,\dots,\down_\pdim)$ and $\up =
  (\up_0,\dots,\up_\pdim)$, let $\vecsp$ be the vector space over the
  rational numbers of all matrices $(\beta_{i,j})$ such that:
  \begin{enumerate}
  \item[a)] $\beta_{i,j}$ is a solution to the system of linear
    equations (\ref{vSpaceEq}),
  \item[b)] $\beta_{i,j}=0$ whenever $j<\down_i$ or $j>\up_i$ (or $i <
    0$ or $i>p$).
  \end{enumerate}
  An element $D = (\beta_{i,j})$ in $\vecsp$ will be called a
  $\emph{diagram}$ and we define the $h$-vector, $h_D(t)$, and
  multiplicity, $e(D)$, of a diagram from the polynomial $S_D(t) =
  \sum_{i,j} (-1)^i \beta_{i,j} t^j$ in the same way as we did for a
  module above. Furthermore, we say that $D$ has codimension $\pdim$.
\end{dfn}

Let $M$ be a Cohen-Macaulay $R$-module of codimension $\pdim$.  The
\emph{maximal and minimal shifts} of degree $i$ of $M$ are defined by
$\up_i(M) = \max\{j\,|\, \beta_{i,j}(M) \neq 0\}$ and $\down_i(M) =
\min\{j\,|\, \beta_{i,j}(M) \neq 0\}$, respectively. It is well known
that when $M$ is Cohen-Macaulay then the sequences $\up(M) =
(\up_0(M),\up_1(M),\dots,\up_\pdim(M))$ and $\down(M) =
(\down_0(M),\down_1(M),\dots,\down_\pdim(M))$ are strictly increasing.
Using the notation introduced in Definition 2.1 we get that the Betti
diagram of $M$, $\beta(M)$, is an element in the vector space
$V_{\down(M),\up(M)}$.

Let $\down$ and $\up$ be two strictly increasing sequences of integers
and denote by $\bettiset$ the set of Betti diagrams in $\vecsp$. As
mentioned in the introduction we have that $\beta(M_1) + \beta(M_2) =
\beta(M_1 \oplus M_2)$ so we see that $\beta(M_1) + \beta(M_2)$ is in
$\bettiset$ if $\beta(M_1)$ and $\beta(M_2)$ are, for any
Cohen-Macaulay $R$-modules $M_1$ and $M_2$. By normalizing the Betti
diagrams in $\bettiset$ in such a way that $\beta_0 = \sum_j
\beta_{0,j} = 1$ we get a subset of $\vecsp$ closed under convex
combinations over the rational numbers. By describing this convex, and
in fact bounded, subset of the vector space $\vecsp$ we get a
description of the possible Betti diagrams up to multiplication by an
integer.

\begin{dfn}
  Let $M$ be a Cohen-Macaulay $R$-module.  The numbers
  $\cbetti_{i,j}(M)=\beta_{i,j}(M)/\beta_0$, where $\beta_0 =
  \sum_j\beta_{0j}(M)$, will be called the \emph{normalized Betti
    numbers} of $M$ and the corresponding diagram $\cbetti(M) =
  \beta(M)/\beta_0$ the \emph{normalized Betti diagram} of $M$.  We
  denote by $\cbettiset$ the set of all normalized Betti diagrams in
  $\vecsp$.
\end{dfn}


\subsection{Pure diagrams}
The main purpose of this section is to give a conjectural description
of the set of normalized Betti diagrams, $\cbettiset$, as the convex
hull of the pure diagrams in $\vecsp$ and we begin by explaining the
notion of a pure diagram.

We say that a $R$-module $M$ has a \emph{pure resolution} of type
$d=(d_0,d_1,\dots, d_\pdim)$ if its minimal free resolution has the
form
$$
0 \to R(-d_\pdim)^{\beta_{\pdim}} \to \dots \to R(-d_0)^{\beta_{0}}
\to M \to 0.
$$
Herzog and K\"{u}hl showed in \cite{Herzog-Kuhl} that the Betti
numbers of a pure resolution of a Cohen-Macaulay algebra are
determined by its type and Huneke and Miller computed in
\cite{Huneke-Miller} the multiplicity of such an algebra, also in
terms of its type. These two results can be extended to Cohen-Macaulay
$R$-modules, see \cite{Migliore-Nagel-Roemer2}, and what follows is a
short explanation of how this can be done.

Given a pure resolution we may try to solve system (\ref{vSpaceEq})
for the unknowns $\beta_0, \beta_1, \dots, \beta_\pdim$.  Before we do
that we will add one equation.  We have
$$
\frac{d^p}{dt^\pdim} S_M(t) = \frac{d^\pdim}{dt^\pdim}
h_M(t)(1-t)^\pdim = \sum_{i=0}^\pdim \binom{\pdim}{i}(-1)^i
\frac{\pdim!(1-t)^{\pdim-i}}{(\pdim-i)!}
\frac{d^{\pdim-i}}{dt^{\pdim-i}}h_M(t)
$$
where the second equality follows by expanding its left-hand side with
the product rule for derivation.  With $t=1$ this equation yields
$\frac{d^p}{dt^\pdim} S_M(1) = (-1)^{p}p!h_M(1) = (-1)^{p}p!e(M)$.
Since, on the other hand, $\frac{d^\pdim}{dt^\pdim} S_M(1) =
\frac{d^\pdim}{dt^\pdim} \sum_{i,j} (-1)^i \beta_{i,j}(M) t^j =
\sum_{i,j} (-1)^i j^{\pdim}\beta_{i,j}$, where we have used
(\ref{vSpaceEq}) to get rid of the the lower powers of $j$, we get
$$
\sum_{i,j} (-1)^i j^{\pdim}\beta_{i,j} = (-1)^\pdim\pdim!e(M).
$$
Adding this equation to the system (\ref{vSpaceEq}) yields
$$
\begin{pmatrix}
  1 & -1 & \dots & (-1)^{\pdim} \\
  d_0 & -d_1 & \dots & (-1)^{\pdim}d_\pdim \\
  \vdots & \vdots & \vdots & \vdots \\
  d_0^{\pdim-1} & -d_1^{\pdim-1} & \dots & (-1)^{\pdim}d_\pdim^{\pdim-1} \\
  d_0^\pdim & -d_1^\pdim & \dots & (-1)^{\pdim}d_\pdim^\pdim
\end{pmatrix}
\begin{pmatrix} \beta_0 \\ \beta_1 \\ \vdots \\ \beta_{\pdim-1} \\
  \beta_{\pdim} \\\end{pmatrix} =
\begin{pmatrix} 0 \\ 0 \\ \vdots \\ 0 \\ (-1)^\pdim \pdim!e(M)
\end{pmatrix}
$$
and by solving this linear system of equations with Cramer's rule and
the Vandermonde determinant we get
$$
\beta_i = (-1)^{i+1} \pdim!e(M) \prod_{k\neq i} \frac{1}{d_k-d_i}.
$$
Hence
$$
e(M) = \beta_0\frac{1}{\pdim!} \prod_{k=1}^\pdim (d_k-d_0)
$$
and
$$
\beta_i = \beta_0(-1)^{i+1} \prod_{\substack{k \neq i \\k \neq 0}}
\frac{d_k-d_0}{d_k-d_i}.
$$

Note that this solution exists for any strictly increasing sequence of
integers $(d_0,d_1,\dots,d_p)$, even if it is not given by the type of
a pure resolution.

\begin{dfn}
  For any strictly increasing sequence of integers $d=(d_0,d_1,\dots,
  d_\pdim)$, the matrix $\pure(d)$ defined by
$$
\pure(d)_{i,j} =
\begin{cases}(-1)^{i+1} \prod_{\substack{k \neq i \\k \neq 0}}
  \frac{d_k-d_0}{d_k-d_i}
  \text{ when $j = d_i$,} \\
  0 \text{ when $j \neq d_i$},
\end{cases}
$$
will be called a \emph{pure diagram}. The partial order on sequences
of integers of length $\pdim+1$, given by $d \leq d^{\prime}$ whenever
$d_i \leq d^{\prime}_i$ for $i=0,1,\dots,\pdim$, defines a partial
order on the pure diagrams, that is, $\pure(d) \leq \pure(d^{\prime})$
whenever $d \leq d^{\prime}$. Furthermore, for any two strictly
increasing sequences, $\down$ and $\up$, of length $\pdim+1$ we denote
by $\poset$ the partially ordered set of all pure diagrams in
$\vecsp$, that is, all $\pure(d)$ such that $\pure(\down) \leq
\pure(d) \leq \pure(\up)$.
\end{dfn}

Note that the multiplicity of a pure diagram, by the above
calculation, is given by
$$
e(\pure(d)) = \frac{1}{\pdim!} \prod_{i=1}^\pdim (d_i-d_0),
$$
and hence that $\pure(d) < \pure(d^\prime)$ implies $e(\pure(d)) <
e(\pure(d^\prime))$ as long as $d_0 = d^{\prime}_0$.

We are now ready to state the conjecture describing the set of
normalized Betti diagrams $\cbettiset$.

\begin{conj}\label{MainConj}
  The Betti diagram of any Cohen-Macaulay $R$-module is a non-negative
  linear combination of pure diagrams and furthermore, any pure
  diagram is a rational multiple of the Betti diagram of some
  Cohen-Macaulay $R$-module. In other words, $\cbettiset$ is equal to
  the convex hull of $\poset$.
\end{conj}

\begin{example} \label{ConjEx}
  Let $A = k[x,y]/(y^2,xy,y^4)$. The Betti diagram of $A$ is then
$$
\beta(A) = \begin{pmatrix} 1&-&- \\ -&2&1 \\ -&-&- \\ -&1&1
\end{pmatrix}
$$
where we have used the convention of writing the entry
$\beta_{i,j}(A)$ in column $i$ and row $j-i$. (The reason for shifting
the rows in this way is simply to save space.)

Conjecture \ref{MainConj} states that any Betti diagram of a
Cohen-Macaulay $R$-module allows a non-negative linear expansion into
pure diagrams and we see that this holds for the Betti diagram
$\beta(A)$ since
$$
\beta(A) =\frac{1}{2}
\begin{pmatrix} 1&-&- \\ -&3&2 \\ -&-&- \\ -&-&- \end{pmatrix}
+\frac{3}{10}
\begin{pmatrix} 1&-&- \\ -&5/3&- \\ -&-&- \\-&-&2/3 \end{pmatrix}
+\frac{1}{5}
\begin{pmatrix} 1&-&- \\ -&-&- \\ -&-&- \\ -&5&4 \end{pmatrix},
$$
that is,
$$
\beta(A) =
\frac{1}{2}\pure(0,2,3)+\frac{3}{10}\pure(0,2,5)+\frac{1}{5}\pure(0,4,5).
$$

The conjecture also states that some integer multiple of any pure
diagram is the Betti diagram of a Cohen-Macaulay $R$-module. We will
now see that this is true for the pure diagrams in this example.  The
pure diagrams $\pure(0,2,3)$ and $\pure(0,4,5)$ are the Betti diagrams
of $k[x,y]/\mathfrak{m}^2$ and $k[x,y]/\mathfrak{m}^4$ where
$\mathfrak{m} = (x,y)$. The pure diagram $\pure(0,2,5)$ have
non-integer entries so it needs to be multiplied by some integer. Let
$I$ and $J$ be the ideals of $k[x,y]$ given by $I = (x^4,x^2y^2,y^4)$
and $J = (x^6,x^3y^3,y^6)$ and let $M$ be the $k[x,y]$-module $M =
I/J$.  Then the Betti diagram of $M$ is, after twisting its degrees by
$4$, $\beta(M(4)) = 3\cdot\pure(0,2,5)$.
\end{example}

We will now show that any maximal chain of pure diagrams is a vector
space basis of $\vecsp$. Hence, we know that any Betti diagram is a
linear combination of pure diagrams while Conjecture~\ref{MainConj}
states that there is a linear combination with non-negative
coefficients.  A natural question now is if the pure diagrams in a
non-negative linear expansion of a Betti diagram always can be chosen
to be a maximal chain.  We make a closer analysis of this question in
Section~\ref{OrderComplex} where we conjecture this to be true
(Conjecture~\ref{ConvexConj}).
 
\begin{prop}
  The elements of any chain of pure diagrams are linearly independent.
\end{prop}

\begin{proof}
  Let $\pure(d_1) < \pure(d_2) < \dots < \pure(d_t)$ be a chain of
  pure diagrams. Then there is an integer $k$ such that $(d_1)_k <
  (d_2)_k \leq (d_3)_k \leq \dots \leq (d_t)_k$ and hence
  $(\pure(d_1))_{k,(d_1)_k} \neq 0$ and $(\pure(d_i))_{k,(d_1)_k} = 0$
  for all $2 \leq i \leq t$. This shows that $\pure(d_1)$ is not in
  the span of the pure diagrams
  $\pure(d_2),\pure(d_3),\dots,\pure(d_t)$. Now we may assume by
  induction that the pure diagrams in the shorter chain $\pure(d_2) <
  \pure(d_3) < \dots < \pure(d_t)$ are linearly independent which
  finishes the proof.
\end{proof}

\begin{prop} 
  For any two strictly increasing sequences of non-negative integers
  $\down = (\down_0, \down_1,\dots, \down_\pdim)$ and $\up =
  (\up_0,\up_1,\dots,\up_\pdim)$, the dimension of $\vecsp$ is
  $1+\sum_{i=0}^\pdim (\up_i-\down_i)$ and any maximal chain in
  $\poset$ is a basis of $\vecsp$.
\end{prop}

\begin{proof}
  If $D$ is a diagram in $\vecsp$ then $D$ has at most
  $\sum_{i=0}^\pdim (\up_i-\down_i+1) = \pdim + 1 + \sum_{i=0}^\pdim
  (\up_i-\down_i)$ non-zero entries.  Let $j_0,j_1,\dots,j_t$ be the
  set of all integers $j$ such that $\down_i \leq j \leq \up_i$ for
  some $i$, and note that $t \geq \pdim$. Then we can write equation
  (\ref{vSpaceEq}) as
$$
\left\{
  \begin{array}{ccccccccc}
    \alpha_{j_0}&+ &\alpha_{j_1}&+&\dots& + &\alpha_{j_t} &=& 0 \\
    j_0 \alpha_{j_0}&+&j_1 \alpha_{j_1}&+&\dots &+&j_t \alpha_{j_t} &=& 0 \\
    &&&&&&& \vdots & \\
    j_0^{p-1} \alpha_{j_0}&+&j_1^{p-1} \alpha_{j_1}&+&\dots& +&j_t^{p-1} \alpha_{j_t} &=& 0 \\
  \end{array}\right.
$$
where $\alpha_{j_k} = \sum_i (-1)^i \beta_{i,j_k}$. This linear system
has maximal rank by the Vandermonde determinant and since $t \geq
\pdim$ this rank is $\pdim$. We get $\dim \vecsp = 1+\sum_{i=0}^p
(\up_i-\down_i)$.

Let $F$ be a maximal chain in $\poset$. The elements of $F$ are
linearly independent by Proposition~2.6. To see that they span
$\vecsp$ it is enough to prove that the length of $F$ is equal to the
dimension of $\vecsp$. Let $\pure(\down) = \pure(d_1) < \pure(d_2) <
\dots < \pure(d_n) = \pure(\up)$ be a maximal chain.  Let the degree
of $d_i$ be defined by $\deg d_i = \sum_k (d_i)_k$.  Since the chain
is maximal we have, for each $i$, that $d_i$ and $d_{i+1}$ differ in
precisely one position and that they differ by one in that position,
that is $\deg d_{i+1} - \deg d_{i} = 1$. We get $\deg \up - \deg \down
= \sum_{i=1}^{n-1} \deg d_{i+1} - \deg d_{i} = n-1$ and hence the
length of the chain is $n = 1 + \deg \up - \deg \down = \dim \vecsp$.
\end{proof}


\subsection{The Multiplicity conjecture}
One of the main motivations of Conjecture~2.4 is the Multiplicity
conjecture of Herzog, Huneke and Srinivasan (see \cite{Huneke-Miller},
\cite{Herzog-Srinivasan}).
Let $A = R/I$ be a Cohen-Macaulay algebra for some ideal $I$ of $R$, and
let $\down = (\down_0, \down_1, \dots, \down_p)$ and 
$\up = (\up_0, \up_1, \dots, \up_p)$ be the minimal and maximal shifts of
a minimal free resolution of $A$,
that is, $\down_i = \min \{j\,|\,\beta_{i,j}(A) \neq 0\}$ and
 $\up_i = \max \{j\,|\,\beta_{i,j}(A) \neq 0\}$.
The Multiplicity conjecture states that
$$
\prod_{i=1}^\pdim \down_i \leq \pdim!e(A) \leq \prod_{i=1}^\pdim \up_i.
$$
Migliore, Nagel and Römer extended this conjecture in \cite{Migliore-Nagel-Roemer2}
by saying that, we have
equality in the above equation, below or above, if and only if the
resolution of $A$ is pure.

We will now show that this extended Multiplicity conjecture holds for
all diagrams in the convex hull of $\poset$ and thus that
Conjecture~2.4 implies the extended Multiplicity conjecture.
This turns out to be an easy consequence of the fact that $\pi > \pi^{\prime}$
implies $e(\pi) > e(\pi^{\prime})$ for any pure diagrams $\pi$ and $\pi^{\prime}$
of the same codimension, whose type agree in degree zero.

\begin{prop}\label{ConjImplMC}
Let $\down = (\down_0, \down_1, \dots, \down_p)$ and 
$\up = (\up_0, \up_1, \dots, \up_p)$ be two strictly increasing sequences of
integers such that $\down \leq \up$.
Assume that $\down_0 = \up_0 = 0$ and that $D \in \vecsp$ is in the 
convex hull of $\poset$.
Then
$$
\prod_{i=1}^\pdim \down_i \leq \pdim!e(D) \leq \prod_{i=1}^\pdim \up_i.
$$
with equality below or above if and only if $D$ is a pure diagram.
\end{prop}

\begin{proof}
We have that
$$
D = \sum_{\pure \in \poset} c_{\pure}\pure
$$ 
for some non-negative rational numbers $c_{\pure}$ such that 
$\sum_{\pure \in \poset}c_{\pure} = 1$.
We get
$$
e(D) =  \sum_{\pure \in \poset} c_\pure e(\pure)
$$
and hence
$$
\min \,\{e(\pure)\,|\,\pure \in \poset\} \leq e(D) \leq
\max \,\{e(\pure)\,|\,\pure \in \poset\}.
$$
Since $\pure < \pure^{\prime}$ implies that 
$e(\pure) < e(\pure^{\prime})$ we get
$e(\pure(\down)) \leq e(D) \leq e(\pure(\up))$ and the inequalities
of the proposition follows
by applying the formula for the multiplicity of a pure diagram to
$e(\pure(\down))$ and $e(\pure(\up))$.
Furthermore, if $D \neq \pure(\down)$ 
then the left inequality is strict which shows that
we have equality below if and only if $D$ is a pure diagram. In the same way we see
that we have equality above if and only if $D = \pure(\up)$.
\end{proof}


\subsection{The order complex of $\poset$} \label{OrderComplex} In
this section we will present a conjecture giving the convex hull of
$\poset$ more structure.

To any partially ordered finite set, \emph{poset} for short, we can
associate a simplicial complex.  The \emph{order complex}
$\Delta(\varPi)$ of the poset $\varPi$ is simply the set of chains in
$\varPi$. The chains of $\varPi$ are thus the \emph{faces} of
$\Delta(\varPi)$ and the maximal chains its \emph{facets}. The order
complex we are interested in is, of course, $\Delta(\poset)$ and we
will begin by describing a geometric realization, in the vector space
$\vecsp$, of this abstract simplicial complex. We then conjecture that
the convex hull of $\poset$ is equal to this geometric realization of
$\Delta(\poset)$

Since the elements in a maximal chain in $\poset$, by Proposition~2.6,
are linearly independent, the convex hull of these elements is a
geometric realization of the abstract simplex corresponding to the
chain. We will now see that the union of these simplices is a
geometric realization of the order complex $\Delta(\poset)$, so we may
think of $\Delta(\poset)$ as sitting in $\vecsp$.

\begin{figure}[htbp]\caption{$\poset$ and $\Delta(\poset)$ when $\down
    = (0,1,3)$ and $\up = (0,3,4)$. }
  \label{fig_1}
  \vspace{10pt} \input{ordercomplex.pstex_t}
\end{figure}

\begin{prop}\label{GeoRealization}
  The simplicial complex $\Delta(\poset)$ has a geometric realization
  given by the union of the simplices spanned by the maximal chains in
  $\poset$.
\end{prop}

\begin{proof}
  Let $\pi_1 < \pi_2 < \dots < \pi_{\dimvecsp}$ and $\pi^{\prime}_1 <
  \pi^{\prime}_2 < \dots < \pi^{\prime}_{\dimvecsp}$ be two maximal
  chains in $\poset$.  It is sufficient to show that any point in the
  intersection of the two simplices spanned by $\pi_1, \pi_2, \dots,
  \pi_{\dimvecsp}$ and $\pi^{\prime}_1,\pi^{\prime}_2, \dots,
  \pi^{\prime}_{\dimvecsp}$ is in a common face, spanned by a subset
  of the common vertices of the two simplices.

  Assume that
$$
\sum_{i=1}^{\dimvecsp} a_i\pi_i = \sum_{i=1}^{\dimvecsp} a_i^\prime
\pi_i^\prime,
$$
for some non-negative rational numbers $a_1,a_2, \dots, a_{\dimvecsp}$
and $a_1^\prime, a_2^\prime, \dots, a_{\dimvecsp}^\prime$.  We will
now prove that $a_i = a_i^\prime$ for all $0 \leq i \leq \dimvecsp$
and furthermore that $a_i = a_i^\prime = 0$ whenever $\pi_i \neq
\pi_i^\prime$.

Assume that this is true for the sub-chains $\pi_1 < \pi_2 < \dots <
\pi_k$ and $\pi^{\prime}_1 < \pi^{\prime}_{2} < \dots <
\pi^{\prime}_{k}$. Then
$$
\sum_{i=k+1}^{\dimvecsp} a_i\pi_i = \sum_{i=k+1}^{\dimvecsp}
a_i^\prime \pi_i^\prime.
$$
In what follows we will refer to the integers in the sequence giving
the type of a pure diagram as its shifts.  Since $\pi_{k+1}$ is the
smallest element in the chain $\pi_{k+1} < \pi_{k+2} < \dots <
\pi_{\dimvecsp}$ it has a shift in some position unique among the
shifts occurring in this position for all the pure diagrams in the
chain. Hence it has a non-zero entry in some position where the entry
of all other pure diagrams in the chain are zero. The same holds for
$\pi_{k+1}^\prime$ and the chain $\pi_{k+1}^\prime < \pi_{k+2}^\prime
< \dots < \pi_{\dimvecsp}^\prime$.  It follows that if $\pi_{k+1} =
\pi_{k+1}^\prime$ then $a_{k+1} = a_{k+1}^\prime$.

If instead $\pi_{k+1} \neq \pi_{k+1}^\prime$, then $\pi_{k+1}$ and
$\pi_{k+1}^\prime$ are incomparable. This means $\pi_{k+1}$ has a
non-zero entry in a position where the entry of all the pure diagrams
in the chain $\pi_{k+1}^\prime < \pi_{k+2}^\prime < \dots <
\pi_{\dimvecsp}^\prime$ are zero and thus $a_{k+1} = 0$. That
$a_{k+1}^\prime = 0$ follows by a symmetric argument.  We have shown
that the assertion holds for the chains $\pi_1 < \pi_2 < \dots <
\pi_{k+1}$ and $\pi_{1}^\prime < \pi_{2}^\prime < \dots <
\pi_{k+1}^\prime$.  The assertion holds obviously for the empty chains
so we can use this as basis for the induction which finishes the
proof.
\end{proof}

Note that Proposition~\ref{GeoRealization} actually shows that any
non-negative linear expansion into a maximal chain is unique, since,
after normalization, it will be contained in precisely one of the
simplices of $\Delta(\poset)$.

\begin{conj} \label{ConvexConj} Any non-negative linear combination of
  pure diagrams is a non-negative linear combination of pure diagrams
  from the same chain. In other words, the convex hull of $\poset$ is
  equal to the geometric realization of the order complex
  $\Delta(\poset)$.
\end{conj}

Conjecture~\ref{ConvexConj} leads to an interesting algorithm for
deciding whether a diagram, $D$, has a non-negative linear expansion
into pure diagrams, or not.  By Conjecture~\ref{ConvexConj} this is
equivalent to deciding if it has an expansion into a chain of pure
diagrams.  Assume that
$$
D = c_1\pure(d_1) + c_2\pure(d_2) + \dots + c_\pdim \pure(d_\pdim)
$$
where $\pure(d_1) < \pure(d_2) < \dots < \pure(d_\pdim)$ is a chain of
pure diagrams and $c_1, c_2, \dots , c_\pdim$ are non-negative
rational coefficients.  Then $d_1$ is the sequence of minimal shifts
of $D$ and $c_1$ is the largest non-negative integer such that all the
entries of $D^\prime = D - c_1\pure(d_1)$ are non-negative. We get
$$
D^\prime = c_2\pure(d_2) + \dots + c_\pdim \pure(d_\pdim)
$$
and continuing in the same way gives us the chain $\pure(d_1) <
\pure(d_2) < \dots < \pure(d_\pdim)$ and the coefficients $c_1, c_2,
\dots , c_\pdim$. If $D$ is not a non-negative linear combination of
pure diagrams from the same chain then, at some point, the minimal
shifts of $D^\prime$ will be a sequence of integers that is not
strictly increasing.

\begin{example}
Consider the Betti diagram
$$
\beta(A) = \begin{pmatrix} 1&-&- \\ -&2&1 \\ -&-&- \\ -&1&1
\end{pmatrix}
$$
of Example~\ref{ConjEx}. We will now apply the algorithm described
above to the diagram $\beta(A)$. Its minimal shifts are $(0,2,3)$
and following the algorithm we see that $c_1 = \frac{1}{2}$ is 
the largest non-negative integer such that all entries of
$\beta(A)-c_1\pure(0,2,3)$ are non-negative.
We get
$$
D' = \beta(A)-\frac{1}{2}\pure(0,2,3)
=\begin{pmatrix} 1/2&-&- \\ -&1/2&- \\ -&-&- \\ -&1&1
\end{pmatrix}
$$
whose minimal shifts are $(0,2,5)$. Removing as much as possible
of $\pure(0,2,5)$ from $D'$ without getting negative entries
yields
$$
D'' = D' - \frac{3}{10}\pure(0,2,5)
= \begin{pmatrix} 1/5&-&- \\ -&-&- \\ -&-&- \\ -&1&4/5
\end{pmatrix}
$$
Since $D'' = \frac{1}{5}\pure(0,4,5)$ we see that
$$
\beta(A) =
\frac{1}{2}\pure(0,2,3)+\frac{3}{10}\pure(0,2,5)+\frac{1}{5}\pure(0,4,5).
$$
as in Example~\ref{ConjEx}.
\end{example}


\subsection{The boundary of  $\Delta(\poset)$ and quasipure resolutions}
The boundary of $\Delta(\poset)$ is given by all
codimension one faces that are contained in precisely one facet.
Each such face, $F_i$, defines a halfspace $H_i$ and the intersection
of these halfspaces defines a polytope
$$
\cap_i H_i.
$$
It is clear that $\cap_i H_i \subseteq \Delta(\poset)$, with equality
if and only if $\Delta(\poset)$ is equal to the convex hull of
$\poset$. Hence to prove Conjecture 2.10 it is enough to prove that
$\cap_i H_i = \Delta(\poset)$.
Some of these halfspaces consists of all diagrams where
a specific entry is non-negative and some do not.

\begin{prop}
Let $\pi = \pure(d_1) < \dots < \pure(d_{i-1}) < 
\pure(d_i) < \pure(d_{i+1}) < \dots < \pure(d_{\dimvecsp})$
be a maximal chain in $\poset$. Let $F$ be the codimension one
face in $\Delta(\poset)$ corresponding to the chain obtained by
removing $\pure(d_i)$ from $\pi$. Then $F$ belongs to the boundary
of $\Delta(\poset)$ if and only if one of the following things hold:
\begin{enumerate}
\item[a)] $i = 1$ or $i = {\dimvecsp}$,
\item[b)] $d_{i-1}$ and $d_{i+1}$ differ in only one position,
\item[c)] $d_{i-1}$ and $d_{i+1}$ differ in two adjacent positions,
say $k$ and $k+1$, and $(d_{i-1})_{k+1} = (d_{i-1})_k + 1$.
\end{enumerate}
Furthermore, the faces obtained in this way from
any two maximal chains containing $\pure(d_{i-1}) < 
\pure(d_{i}) < \pure(d_{i+1})$ corresponds to
the same halfspace $H$, and in the case of a) and b) 
this halfspace consists of all diagrams in $\vecsp$
having a specific entry non-negative.
\end{prop}

\begin{proof}
The face $F$ belongs to the boundary of $\Delta(\poset)$
precisely when the facet $\pi$ is the only facet containing
$F$. This happens precisely when $\pure(d_i)$ is the only
pure diagram that extends the chain $F$ and we need to
show that this happens only in the situations described
in a), b) and c).

Since $\pure(d_1)$ and $\pure(d_{\dimvecsp})$ are the unique
minimal and maximal elements, respectively, in $\poset$,
removing one of them from $\pi$ gives a face belonging
to the boundary. This shows that $F$ belongs to the boundary
in the case of a). 
Assume that $i \neq 1$, $i \neq \dimvecsp$ and that
$$
d_{i-1} = (a_1, \dots,  a_{\pdim}).
$$
We have that $d_{i-1}$ and $d_{i+1}$ differ in at most two
positions. If they differ in only one position, say $k$, then
$$
d_{i+1} = (a_1, \dots, a_k+2, \dots, a_{\pdim})
$$
and there is only one element between  $d_{i-1}$ and $d_{i+1}$, namely
$$
(a_1, \dots, a_k+1, \dots, a_{\pdim}).
$$
If $d_{i-1}$ and $d_{i+1}$ differ in two positions, say $k$ and $l$,
then
$$
d_{i+1} = (a_1, \dots, a_k+1, \dots, a_l+1, \dots, a_{\pdim})
$$
and whenever 
$a_l \neq a_k+1$, there are 
two elements between $d_{i-1}$ and $d_{i+1}$, namely
$$
(a_1, \dots, a_k+1, \dots, a_l, \dots, a_{\pdim})
$$
and
$$
(a_1, \dots, a_k, \dots, a_l+1, \dots, a_{\pdim}).
$$
If instead $a_l = a_k + 1$ then there is only on element between
$d_{i-1}$ and $d_{i+1}$ and we must have 
$k = l + 1$ since the sequence is strictly increasing.
Hence a), b) and c) are the only cases where we get a face
belonging to the boundary of $\Delta(\poset)$.

Let $\pi$ and $\pi^\prime$ be two maximal chains in $\poset$ containing
$\pure(d_{i-1}) < \pure(d_{i}) < \pure(d_{i+1})$. We know from
Proposition~2.6 that the subchains of $\pi$ and $\pi^\prime$
of all pure diagrams between $\pure(\down)$ and $\pure(d_{i-1})$
span the same subspace of $\vecsp$, namely $V_{\down,d_{i-1}}$
and the same thing holds for the subchains of all pure
diagrams between  $\pure(d_{i+1})$ and $\pure(\up)$. Hence,
when we remove $\pure(d_{i})$ from $\pi$ and $\pi^\prime$ the
remaining pure diagrams in the two chains span the same
subspace of $\vecsp$. Hence the two faces of $\Delta(\poset)$
defined by removing $\pure(d_{i})$ from $\pi$ and $\pi^\prime$ define the
same halfspace.

In the case of a) and b), $d_{i}$ has an entry in a position, say $k$,
unique among all entries in position $k$ of the other shifts
in the chain. This means that if $D \in \vecsp$ then, since the chain
$\pi$ is a basis, $D$ is a linear combination of the
pure diagrams in $\pi$ and the coefficient of $\pure(d_i)$
is zero precisely when $D_{k,(d_i)_k}$ is. This means that
the subspace of $\vecsp$ spanned by the chain obtained by removing $\pure(d_i)$
from the chain $\pi$ consists of the diagrams $D$ such that $D_{k,(d_i)_k} = 0$.
Hence the corresponding halfspace consists of all diagrams in which the entry
in position $(k,(d_i)_k)$ is zero.
\end{proof}

If a resolution has its Betti diagram in $\bettiset$, where
$\up_{i-1} \leq \down_i$ for all $0 \leq i \leq \pdim$, it is called \emph{quasipure}.
Herzog and Srinivasan showed in \cite{Herzog-Srinivasan}
that the multiplicity conjecture holds for quasipure resolutions.
As a corollary to Proportion~2.11 we get that part of Conjecture~2.4
holds for a \emph{strictly quasipure resolutions}, that is, resolutions
with minimal and maximal shifts satisfying
$\up_{i-1} < \down_i$ for all $0 \leq i \leq \pdim$.

\begin{cor}
The Betti diagram of any Cohen-Macaulay module with strictly quasipure
resolution is a non-negative linear combination of pure diagrams.
In other words, if $\up_{i-1} < \down_i$ for all $0 \leq i \leq \pdim$ then
$\cbettiset \subseteq \Delta(\poset)$. It is also true, in this case, that
 $\Delta(\poset)$ is equal to the convex hull of $\poset$.
\end{cor}

\begin{proof}
The point is that from the hypothesis on $\up$ and $\down$ it follows
that the codimension one faces on the boundary of $\Delta(\poset)$
are all of the type described in a) and b) in Proposition~2.11.
This means that if $H_1,H_2,\dots,H_t$ are the halfspaces
corresponding to the faces on the boundary, then each $H_i$ consists of
all diagram in $\vecsp$ with a specific entry non-negative.
Hence, every diagram in $\vecsp$ with all entries non-negative
is in $\cap_i H_i$. Since each diagram in $\poset$
 have  non-negative entries 
we see that $\poset \subseteq \cap_i H_i$ and, as mentioned in the beginning of
this section, this implies that  the convex hull of $\poset$ is equal to $\Delta(\poset)$.
Finally, since the diagrams in $\bettiset$ have non-negative entries
we get $\cbettiset \subseteq \cap_i H_i \subseteq \Delta(\poset)$.
\end{proof}


\subsection{Maximal Betti diagrams}
Let $M = F/N$ be a Cohen-Macaulay $R$-module, where $F$ is a free
$R$-module and $N$ a submodule of $F$.  Then there is a lexicographic
submodule $L$ of $F$ such that $F/L$ and $M$ have the same $h$-vector
(see \cite{Macaulay} and \cite{Hulett-2}).  The Betti diagram of $F/L$
is completely determined by the $h$-vector, $h(t)$, of $M$ and of the
free module $F$ and we denote this diagram by $\MaxLex_F(h(t))$.  The
Betti numbers of $M$ are smaller than the Betti numbers of $F/L$, and
we even have that $\beta(M)$ is obtained from $\MaxLex_F(h(t))$ by a
sequence of consecutive cancellations (see \cite{Peeva}).  A
consecutive cancellation is defined as follows.  Let $k$ and $l$ be
two integers and define a diagram $C^{k,l}$ by
$$
(C^{k,l})_{i,j} = \begin{cases} 1 &\text{ when $(i,j) = (k,l)$ or
    $(i,j)=(k+1,l)$}, \\ 0 &\text{ otherwise}.\end{cases}
$$
A consecutive cancellation in position $(k,l)$ of a diagram $D$ is a
diagram
$$
D^{\prime} = D - bC^{k,l},
$$
where $b$ is a non-negative rational number, and we assume that
$D^\prime$ have no negative entries.  Note that a consecutive
cancellation does not change the polynomial $S_{D}(t) = \sum_{i,j}
(-1)^i D_{i,j} t^j$, that is, $S_{D}(t) = S_{D^{\prime}}(t)$, and
hence not the $h$-vector or multiplicity of $D$. In short, for any
Cohen-Macaulay $R$-module $M$ with $h$-vector $h(t)$ we have that
\begin{equation} \label{MaxLexBettiEq} \beta(M) = \MaxLex_F(h(t)) -
  \sum_{i,j} b_{i,j} C^{i,j}
\end{equation}
for some non-negative integers $b_{i,j}$.

Let $F$ be a free $R$-module with basis $e_1,e_2,\dots,e_t$ and assume
that the basis of $F$ is ordered by $e_1 < e_2 < \dots < e_t$ where
$\deg e_1 \leq \deg e_2 \leq \dots \leq \deg e_t$.  Let
$H:\,\mathbb{Z} \to \mathbb{Z}$ be the Hilbert function of some
submodule of $F$. Then, as mentioned above, there is a lexicographic
submodule $L$ with $H$ as Hilbert function and we will now describe
$L$. The monomials of $F$ are the elements on the form $ue_i$ for some
monomial $u$ of $R$ and basis element $e_i$ of $F$.  The lexicographic
order on monomials of $F$ is given by $ue_i > ve_j$ when $i<j$ or when
$i=j$ and $u>v$ in the lexicographic order on monomials of $R$.  The
lexicographic submodule $L$ is generated in each degree $d$ by the
$H(d)$ largest monomials in $F_d$.

For any degree $d$, there are unique integers $k$ and $q$ such that $0
\leq q < H(R,d-\deg e_k)$ and
\begin{equation} \label{LexEq} H(d) = \sum_{i=1}^{k-1} H(R,d-\deg e_i)
  + q.
\end{equation}
For any $F$, $H$ and $d$ denote these numbers by $k(F,H,d)=k$ and
$q(F,H,d) = q$ and furthermore let $g(F,H,d) = \deg e_k$.  This means
that $L$ is generated in degree $d$ by
$$
\sum_{i=1}^{k-1} \mathfrak{m}^{d-\deg e_i}e_i + Ie_{k}
$$
where $k = k(F,H,d)$ and $I$ is the ideal generated by the $q(F,H,d)$
largest monomials of degree $d-g(F,H,d)$ in $R$.

In what follows, let $e_{i,j}$, where $i=1,2,\dots,t$, $j =
1,2,\dots,s$ and $\deg e_{i,j} = \deg e_i$, be a basis of $F^s$
ordered by $e_{i,j} > e_{i',j'}$ when $i < i'$ or when $i = i'$ and $j
> j'$, for any integer $s$.

\begin{lemma} \label{LexLemma} Let $N$ be a submodule of $F$ such that
  the module $M = F/N$ is artinian. Then there is an integer $s$ and a
  lexicographic submodule $L$ of $F^s$ such that $L$ and $N^s$ have
  the same $h$-vector, and $L$ is on the form
$$
L = \sum_{j=1}^s\sum_{i=1}^t \mathfrak{m}^{c_{i,j}}e_{i.j}
$$
where $c_{i,j}$ are non-negative integers.
\end{lemma}

\begin{proof}
  Let $H$ be the Hilbert function of $M$.  To prove the lemma it is
  enough to show that there is a positive integer $s$ such that
  $q(F^s,sH,d) = 0$ for all non-negative integers $d$.  Note that in
  (\ref{LexEq}), $H(d) = H(F,d)$ implies $q=0$, and hence that $H(d)=
  H(F,d)$ implies $q(F,H,d) = 0$. Since $M$ is artinian, $H(d) \neq
  H(F,d)$ only for a finite number of integers and hence the same
  holds for $q(F,H,d)$. This means that there exists a positive
  integer $s$ such that
  \begin{equation} \label{sEq} s\frac{q(F,H,d)}{H(S,d-g(F,H,d))}
  \end{equation}
  is an integer for all $d$. Let $s$ be such an integer.

  Fix a degree $d$ and let $k = k(F,H,d)$ and $q = q(F,H,d)$ then we
  get by multiplying both sides of (\ref{LexEq}) by $s$,
$$
sH(d) = s\sum_{i=1}^{k-1} H(R,d-\deg e_i) + sq =
\sum_{j=1}^s\sum_{i=1}^{k-1} H(R,d-\deg e_{i,j}) + sq
$$
and since by (\ref{sEq}), $sq = s'H(R,d-\deg e_k)$ for some $0 \leq s'
< s$ we get
$$
sH(d) = \sum_{j=1}^s\sum_{i=1}^{k-1} H(R,d-\deg e_{i,j}) +
\sum_{i=j}^{s'}H(R,d-\deg e_{k,j}).
$$
From the above equation it follows that $q(F^s,sH,d) = 0$ for all
non-negative integers $d$. This means that $L$ is generated in degree
$d$ by
$$ \sum_{j=1}^s\sum_{i=1}^{k-1} \mathfrak{m}^{d-\deg e_{i,j}}e_{i,j} +
\sum_{i=j}^{s'} \mathfrak{m}^{d-\deg e_{k,j}}e_{i,j}.
$$
and since this holds for any degree $d$ the lemma follows.
\end{proof}

\begin{prop}\label{NormMaxBettiProp}
  Let $M = F/N$ be an artinian $R$-module of codimension $p$ with
  $h$-vector $h(t)$.  Then there exist non-negative rational numbers
  $a_{i,j}$ and an integer $s$ such that the diagram
$$
D = \sum_{i,j} a_{i,j} \beta(R/\mathfrak{m}^{j+1}(-\deg e_i))
$$
have $h$-vector $h(t)$ and
\begin{equation*}
  \frac{1}{s}\MaxLex_{F^s}(sh(t)) = D.
\end{equation*}
As a consequence we get that $\beta(M)$ is obtained from $D$ by a
sequence of consecutive cancellations, that is,
$$
\beta(M) = D - \sum_{i,j}b_{i,j}C^{i,j}
$$
for some non-negative rational numbers $b_{i,j}$.
\end{prop}

\begin{proof}
  By Lemma~\ref{LexLemma} there is an integer $s$ and a lexicographic
  submodule $L$ of $F^s$ such that $F^s/L$ have $h$-vector $sh(t)$ and
  $L$ is on the form
$$
L = \sum_{j=1}^s\sum_{i=1}^t \mathfrak{m}^{c_{i,j}}e_{i,j}
$$ 
where $c_{i,j}$ are non-negative integers.  We get
$$
\MaxLex_{F^s}(sh(t)) = \beta(F^s/L)
$$
and since then
$$
F^s/L =
\bigoplus_{j=1}^s\bigoplus_{i=1}^tR/\mathfrak{m}^{c_{i,j}}(-\deg
e_{i,j})
$$
we see that
$$
\MaxLex_{F^s}(sh(t)) = \sum_{j=1}^s\sum_{i=1}^t
R/\mathfrak{m}^{c_{i,j}}(-\deg e_{i,j})
$$
Collecting all the terms in the sum above where the $c_{i,j}$ and
$\deg e_{i,j}$ are the same, yields
$$
\MaxLex_{F^s}(sh(t)) = \sum_{i,j}a'_{i,j}
\beta(R/\mathfrak{m}^{j+1}(-\deg e_{i}))
$$
for some non-negative integers $a'_{i,j}$.  Now, let $a_{i,j} =
a'_{i,j}/s$ and $D = \MaxLex_{F^s}(sh(t))/s$.

Since the Hilbert function of $M^s$ is $sh(t)$ we get that
$\beta(M^s)$ is obtained from $\MaxLex_{F^s}(sh(t)) = sD$ by a
sequence of consecutive cancellations, that is,
$$
\beta(M^s) = sD - \sum_{i,j}b'_{i,j}C^{i,j}
$$
for some non-negative integers $b_{i,j}$. Now $\beta(M^s) = s\beta(M)$
so dividing both sides of the equation above with $s$ shows that
$\beta(M)$ is obtained from $D$ by a sequence of consecutive
cancellations.
\end{proof}

\begin{cor} \label{PureMaxBetti} Let M be a Cohen-Macaulay module of
  codimension $p$ generated in degrees $g_1,g_2,\dots,g_t$. Then the
  Betti diagram of $M$ is given by
$$
\beta(M) = \sum_{i,j}a_{i,j}\pure(g_i,j+1,j+2,\dots,j+p)
-\sum_{i,j}b_{i,j}C^{i,j},
$$
for some non-negative rational numbers $a_{i,j}$ and $b_{i,j}$.
\end{cor}

\begin{proof}
  By artinian reduction we may assume that $M$ is artinian. We may
  also assume that $M = F/N$ where $F$ is a free module generated by
  $e_1,e_2,\dots,e_t$ such that $\deg e_i = g_i$ for $i=1,2,\dots,t$.
  By proposition~\ref{NormMaxBettiProp} we get
$$
\beta(M) = D -\sum_{i,j}b_{i,j}C^{i,j},
$$
where
$$
D = \sum_{i,j} a'_{i,j} \beta(R/\mathfrak{m}^{j+1}(-g_i)).
$$
Since
$$
\beta(R/\mathfrak{m}^{j+1-g_i}(-g_i)) = \pure(g_i,j+1, j+2,\dots, j+p)
$$
the corollary follows if we let $a_{i,j} = a'_{i,j-g_i}$.
\end{proof}


\section{Codimension one and two}
We will prove that Conjecture~\ref{MainConj} and
Conjecture~\ref{ConvexConj} holds for Cohen-Macaulay $R$-modules of
codimension one, with generators of any degrees, and for modules of
codimension two, generated in a single degree.

The existence of $R$-modules with pure resolutions, of codimension one
and two, of any type follows from the existence of compressed level
modules.

\begin{prop} \label{Compressed} If $d = (d_0,d_1,\dots,d_p)$ is a
  strictly increasing sequence of integers such that $d_i = d_1+i-1$
  for all $1 \leq i \leq p-1$, then there is a Cohen-Macaulay
  $R$-module with a pure minimal free resolution of type $d$.
\end{prop}

\begin{proof}
  Let $S = k[x_1,x_2, \dots, x_p]$ and $s_i = \dim_k S_i$. Let $t$ and
  $t^{\prime}$ be two positive integers and define
$$
h_i = \min \{ts_i, t^\prime s_{c-i}\}
$$
for each $0 \leq i \leq c$, where $c = d_p - p$.  In earlier work, the
first author showed the existence of compressed level modules in
\cite[Proposition 3.5]{Boij}, and hence that there is an artinian
level $S$-module, $M$, with Hilbert function $i \mapsto h_i$ for any
$t$ and $t^\prime$.

Let $t = s_{c-d_1+1}$ and $t^\prime = s_{d_1-1}$ and let $M$ be a
compressed level $S$-module with Hilbert function $i \mapsto h_i$. Now
it follows from \cite[Proposition 4.1]{Boij} that a minimal free
resolution of $M$ has shifts given by $d = (d_0,d_1,\dots,d_p)$.
\end{proof}

\begin{rmk}
  In codimension one and two there are more explicit ways to construct
  modules with pure resolutions of any type, and these construction
  even gives multigraded modules.  In codimension one everything is
  very simple.  If $A = k[x]/x^{(d_1-d_0)}$ then $\beta(A(-d_0)) =
  \pure(d_0,d_1)$.  In codimension two we can do the following. For
  any strictly increasing sequence of integers, $(d_0,d_1,d_2)$, let
  $a = (d_1-d_0)(d_2-d_1-1)$ and $b = (d_1-d_0)(d_2-d_1)$.  Define two
  ideals of $k[x,y]$ by $I =
  (x^{a-i(d_1-d_0)}y^{i(d_1-d_0)})_{i=0}^{d_2-d_1-1}$ and $J =
  (x^{b-i(d_2-d_1)}y^{i(d_2-d_1)})_{i=0}^{d_1-d_0}$ and let $M$ be the
  multigraded $k[x,y]$-module $M = I/J$. Then $\beta(M(a-d_0)) =
  (d_2-d_1)\pure(d_0,d_1,d_2)$.
\end{rmk}

With the existence part out of the way we turn to the problem of
expanding Betti diagrams into pure diagrams. In codimension one, the
set of possible Betti diagrams is easy to describe and it is then easy
to see that any of them is a non-negative linear combination of pure
diagrams. It turns out that the codimension two case then follows from
the codimension one case. To see this we will need the following
lemma.

For any sequence of integers
$d=(d_0,d_1,\dots,d_{k-1},d_k,d_{k+1},\dots,d_p)$ denote by
$\tau_k(d)$ the sequence $(d_0,d_1,\dots,d_{k-1},d_{k+1},\dots,d_p)$.

\begin{lemma} \label{IsoLemma} Let $k$ be an integer such that and
  $\up_k = \down_k$.  Then there is an isomorphism of vector spaces
$$
\phi_k:\, \vecsp \to V_{\tau_k(\down),\tau_k(\up)}
$$
given by
$$ 
\phi_k(D)_{i,j} = |\up_k-j| \cdot
\begin{cases} D_{i,j} &\text{ when $i<k$,} \\
  D_{i+1,j}& \text{ otherwise,} \end{cases}
$$
for any diagram $D \in \vecsp$.
\end{lemma}

\begin{proof}
  First we need to prove that $\phi_k(D)$ is a diagram of projective
  dimension $\pdim -1$ and by definition this is equivalent to
  $\frac{d^i}{dt^i} S_{\phi_k(D)}(1) = 0$ for $i = 0,1,\dots,p-2$.
  The key observation here is that
$$
S_{\phi_k(D)}(t) = \up_k S_D(t)-\frac{d}{dt} S_D(t).
$$
In fact, since $S_D(t) = \sum_{i,j} (-1)^i D_{i,j} t^j$ and
$\frac{d}{dt} S_D(t) = \sum_{i,j} (-1)^i j D_{i,j} t^j$ we get
$$
\up_k S_D(t)-\frac{d}{dt} S_D(t) = \sum_{i,j} (-1)^i (\up_k -
j)D_{i,j} t^j.
$$
Now
$$
(\up_k-j)D_{i,j} =
\begin{cases} |\up_k-j|D_{i,j} &\text{ when $i<k$,} \\
  0 &\text{ when $i=k$,} \\
  -|\up_k-j|D_{i,j}& \text{ when $i>k$,} \end{cases}
$$
simply because $D_{i,j} = 0$ whenever $j<\down_i$ or $j>\up_i$, and
since by assumption $\up_k = \down_k$.  Finally, we get
\begin{multline*}
  \sum_j \sum_{i=0}^p (-1)^i (\up_k - j)D_{i,j} t^j \\= \sum_j \left(
    \sum_{i = 0}^{k-1} (-1)^i |\up_k-j|D_{i,j}t^j+
    \sum_{i = k+1}^{p} (-1)^{i+1} |\up_k-j|D_{i,j}t^j \right)\\
  = \sum_j \sum_{i=0}^{p-1}(-1)^i \phi_k(D)_{i,j} t^j =
  S_{\phi_k(D)}(t).
\end{multline*}

Since the pure diagrams span $V_{\tau_k(\down),\tau_k(\up)}$ and any
pure diagram, clearly, is in the image of $\phi_k$ we see that
$\phi_k$ is surjective. That $\phi_k$ is an isomorphism follows, for
example, from the fact that $\dim \vecsp = 1+\sum_{i=0}^p
\up_i-\down_i = 1+\sum_{i\neq k} \up_i-\down_i = \dim
V_{\tau_k(\down),\tau_k(\up)}$.
\end{proof}

\begin{thm} \label{CodimOneTwo} Let $M$ be a Cohen-Macaulay $R$-module
  of codimension one, generated in any degrees, or of codimension two,
  generated in a single degree.  Then the Betti diagram of $M$ is a
  non-negative linear combination of pure diagrams from the same chain
  and furthermore, any pure diagram of codimension one or two is a rational
  multiple of the Betti diagram of some Cohen-Macaulay
  $R$-module.  In other words, when $\pdim=1$ or when $\pdim = 2$ and
  $\down_0=\up_0$ then
$$
\cbettiset = \Delta(\poset).
$$
\end{thm}

\begin{proof}
  The existence of $R$-modules of codimension one or two, with pure
  resolutions of any type was proven in Proposition~\ref{Compressed}.

  To prove the rest of the theorem we start with the codimension one
  case.  Let $D$ be a diagram of codimension one. Since for any shifts
  $(d_0,d_1)$ the pure diagram $\pure(d_0,d_1)$ has as its non-zero
  entries $\pure(d_0,d_1)_{0,d_0}=1$ and $\pure(d_0,d_1)_{1,d_1}=1$ we
  see that $D$ is a non-negative linear combination of pure diagrams
  if and only if
  \begin{equation} \label{eqPDimOne_1} D_{i,j} \geq 0
  \end{equation}
  for all integers $i$ and $j$, and
  \begin{equation} \label{eqPDimOne_2} \sum_{i \leq l} D_{0,i} \leq
    \sum_{i \leq l} D_{1,i+1}
  \end{equation}
  for each integer $l$. Furthermore, if $D$ satisfies
  (\ref{eqPDimOne_1}) and (\ref{eqPDimOne_2}) it will be a
  non-negative linear combination of pure diagrams from the same
  chain. To see this let $(d_0,d_1)$ be the minimal shifts of $D$ and
  let $c$ be the largest rational number such that all entries of
  $D^{\prime} = D - c\pure(d_0,d_1)$ are non-negative. Then
  $D^{\prime}$ satisfy (\ref{eqPDimOne_1}) and (\ref{eqPDimOne_2}) and
  the minimal shifts of $D^{\prime}$ are larger than $(d_0,d_1)$.  By
  induction on the minimal shifts we may assume that $D^{\prime}$, and
  hence $D$, is a non-negative linear combination of pure diagrams
  from the same chain.

  Let $M$ be a Cohen-Macaulay $R$-module of codimension one. By
  Corollary~\ref{PureMaxBetti} the Betti diagram of $M$ can be
  obtained from a non-negative linear combination of pure diagrams by
  a sequence of consecutive cancellations. Since a consecutive
  cancellation does not change the properties (\ref{eqPDimOne_1}) and
  (\ref{eqPDimOne_2}) already enjoyed by any non-negative linear
  combination, we see that the Betti diagram of $M$ also satisfies
  (\ref{eqPDimOne_1}) and (\ref{eqPDimOne_2}) and hence is a
  non-negative linear combination of pure diagrams from the same
  chain.
 
  We now turn to the codimension two case. Let $M$ be a Cohen-Macaulay
  $R$-module of codimension two, generated in a single degree. We may,
  without loss of generality, assume that $M$ is generated in degree
  zero.  By Corollary~\ref{PureMaxBetti} we then have $\beta(M) =
  \sum_j a_j \pure(0,j+1,i+2) - \sum_j b_j C^{1,j}$.  Using the
  isomorphism $\phi_0$ from Lemma~\ref{IsoLemma} and the fact that
  $\phi_0(\pure(0,j+1,j+2)) = a\,\pure(j+1,j+2)$, for some rational
  number $a$, and $\phi_0(C^{1,j}) = jC^{0,j}$, we get
  $\phi_0(\beta(M)) = \sum_j a_j^{\prime} \pure(j+1,j+2) - \sum_j b_jj
  C^{0,j}$.  We see that $\phi_0(\beta(M))$ is obtained from a
  non-negative linear combination of pure diagrams by a sequence of
  consecutive cancellations, and since $\phi_0(\beta(M))$ is of
  codimension one, it is, by the same argument as above, a
  non-negative linear combination of pure diagrams from the same
  chain. Finally, $\phi_0$ is a linear map that preserves chains of
  pure diagrams and the same thing holds for its inverse
  $\phi_0^{-1}$.  Since $\beta(M) = \phi_0^{-1}(\phi(\beta(M)))$ we
  get that $\beta(M)$ is a non-negative linear combination of pure
  diagrams from the same chain.
\end{proof}

As a corollary we get, using Proposition~\ref{ConjImplMC}, that the
Multiplicity conjecture holds for Cohen-Macaulay $R$-modules of
codimension two, generated in a single degree.  In
\cite{Migliore-Nagel-Roemer2} the upper bound of the conjecture was
shown to hold for codimension two modules with generators of any
degrees.

\begin{cor}
  If $M$ is a Cohen-Macaulay $R$-module of codimension two, generated
  in a single degree, with minimal and maximal shifts given by $(d_0,
  \down_1, \down_2)$ and $(d_0, \up_1, \up_2)$ respectively, then
$$
\beta_0(M)\frac{(\down_1-d_0)(\down_2-d_0)}{2} \leq e(M) \leq
\beta_0(M)\frac{(\up_1-d_0)(\up_2-d_0)}{2}.
$$
\end{cor}

\section{Higher codimension}
In this section we will show that the Betti diagram of a
Cohen-Macaulay algebra $A$, that is, $A = R/I$ for some ideal $I$ of
$R$, is a non-negative linear combination of pure diagrams when $A$ is
a complete intersection, of any codimension, or a Gorenstein algebra
of codimension three.  The proofs of these two cases have some steps
in common and we will start by explaining these common steps.

Both proofs work by breaking down the Betti diagram into a sum of two
Betti diagrams whose codimension is the codimension of the original
diagram minus one.  Let $E$ and $F$ be two diagrams of codimension
$p-1$ such that $e(E) = e(F)$. Then we claim that $D$, defined by
$$
D_{i,j} = E_{i,j} + F_{i-1,j},
$$
is a diagram of codimension $p$. (Note that $D$ is the sum of the
matrix $E$ and the matrix obtained from $F$ by shifting its columns
one step to the right.)  To prove that $D$ is a diagram of codimension
$p$ we need to prove that the polynomial $S_D(t)$ is divisible by
$(1-t)^p$.  We have
$$
S_D(t) = \sum_{i,j} (-1)^i D_{i,j} t^j = \sum_{i,j} (-1)^i E_{i,j} t^j
+ \sum_{i,j} (-1)^i F_{i-1,j} t^j = S_E(t) - S_F(t).
$$ 
Since $E$ and $F$ are diagrams of codimension $p-1$ we know that
$S_E(t)$ and $S_F(t)$ are divisible by $(1-t)^{p-1}$ and by definition
$h_E(t) = S_E(t)/(1-t)^{p-1}$ and $h_F(t) = S_F(t)/(1-t)^{p-1}$.
Hence $S_D(t)/(1-t)^{p-1} = h_E(t) - h_F(t)$ and it remains to prove
that $t=1$ is a root of $h_E(t) - h_F(t)$. This follows from the fact
that by definition $e(E) = h_E(1)$ and $e(F) = h_F(1)$ and by
assumption $e(E) = e(F)$.

Now assume that $E$ and $F$ are non-negative linear combinations of
pure diagrams.  We will show that if the pure diagrams in these two
non-negative linear combinations fit together in a certain way, then
the diagram $D$ is a non-negative linear combination of pure diagrams.
How two pure diagrams should fit together is described in the
following lemma.

\begin{lemma} \label{stepUp} Let $d = (d_0,d_1,\dots,d_p)$ and
  $d^{\prime} = (d^{\prime}_0,d^{\prime}_1,\dots,d^{\prime}_p)$ be two
  strictly increasing sequences of integers such that $d_i <
  d^{\prime}_i$ for $i=0,1,\dots,p$ and define $D$ by
$$
D_{i,j} = \pure(d)_{i,j} +
\frac{e(\pure(d))}{e(\pure(d^{\prime}))}\pure(d^{\prime})_{i-1,j}.
$$
Then $D$ is a non-negative linear combination of pure diagrams of
codimension $p+1$
\end{lemma}

\begin{proof}
  We may assume, without loss of generality, that $d_0 = 0$.  The
  proof is by induction on the codimension $p$. If $p = 1$, then $D$
  is simply a pure diagram of codimension two with shifts
  $(d_0,d^{\prime}_0)$, so the assertion of the lemma is true in this
  case. Assume now that $p > 1$.  Let $\pi =
  \pure(d_0,d_1,\dots,d_p,d^{\prime}_p)$ and
$$
E = D - \frac{(d^{\prime}_p-d_p)}{d^{\prime}_p}\pi.
$$
Furthermore, let $\phi$ be the the map from Lemma~\ref{IsoLemma}
defined by $\phi(E)_{i,j} = (d^{\prime}_p - j)E_{i,j}$.  We will now
show that
$$
\phi(E)_{i,j} = d_p\pure(\tau(d))_{i,j} +
\frac{e(\pure(d))}{e(\pure(d^{\prime}))}d^{\prime}_p\pure(\tau(d^\prime))_{i-1,j},
$$ 
where $\tau(d) = (d_0,d_1,\dots,d_{p-1})$ and $\tau(d^{\prime}) =
(d^{\prime}_0,d^{\prime}_1,\dots,d^{\prime}_{p-1})$.  Denote by $\rho$
the matrix defined by $\rho_{i,j} = \pure(d^{\prime})_{i-1,j}$.  Then
we have
$$
\phi(\rho)_{i,j} =
(d^{\prime}_p-j)\rho_{i,j}=(d_p^{\prime}-j)\pure(d^{\prime})_{i-1,j}.
$$
Since $\pure(d^{\prime})_{i,d^{\prime}_i} = (-1)^i \prod_{j\neq i}
\frac{d^{\prime}_i - d^{\prime}_0}{d^{\prime}_i-d^{\prime}_j}$ we see
that $(d_p^{\prime}-j)\pure(d^{\prime})_{i-1,j} =
(d^{\prime}_p-d^{\prime}_0)\pure(\tau(d^{\prime}))_{i-1,j}$ and hence
that
$$
\phi(\rho)_{i,j} =
(d^{\prime}_p-d^{\prime}_0)\pure(\tau(d^{\prime}))_{i-1,j}.
$$
In the same way we see that
$$
\phi(\pi)_{i,j} = d^{\prime}_p\pure(d)_{i,j}.
$$
Hence,
\begin{multline*}
  \phi(\pure(d))_{k,d_k} -
  \frac{(d^{\prime}_p-d_p)}{d^{\prime}_p}\phi(\pi)_{k,d_k} =
  (d^{\prime}_p - d_k)\pure(d)_{k,d_k} - (d^{\prime}_p-d_p)\pure(d)_{k,d_k} = \\
  (d_p-d_k)\pure(d)_{k,d_k} = d_p\pure(\tau(d))_{k,d_k}.
\end{multline*}
and we get
$$
\phi(E)_{i,j} = d_p\pure(\tau(d))_{i,j} +
\frac{e(\pure(d))}{e(\pure(d^{\prime}))}(d^{\prime}_p-d^{\prime}_0)\pure(\tau(d^\prime))_{i-1,j}.
$$
Note that
$\frac{e(\pure(d))(d^{\prime}_p-d^{\prime}_0)}{e(\pure(d^{\prime}))d_p}
= \frac{e(\pure(\tau(d)))}{e(\pure(\tau(d^{\prime})))}$ and hence that
$\phi(E)$ is a sum of two pure diagrams in the same way $D$ was, but
now the pure diagrams have codimension $p-1$.  By induction we may
assume that $\phi(E)$ is a non-negative linear combination of pure
diagrams of codimension $p$. By using the inverse of $\phi$ we see
that $E$ is a non-negative linear combination of pure diagrams and
hence that $D$ is, since $D = E + (d^{\prime}_p-d_p)/d^{\prime}_p \pi$
and $d^{\prime}_p-d_p > 0$.
\end{proof}

\begin{example}
  Let $d = (0,2,3)$ and $d' = (1,3,5)$. Then
$$
\pure(d) = \begin{pmatrix} 1&-&-\\-&3&2 \end{pmatrix} \text{ and }
\pure(d') = \begin{pmatrix} -&-&-\\ 1&-&- \\ -&2&- \\-&-&1
\end{pmatrix}
$$
and for the multiplicities we have $e(\pure(d)) = 3$ and $e(\pure(d'))
= 4$.  The diagram
$$
D_{i,j} = \pure(d)_{i,j} +
\frac{e(\pure(d))}{e(\pure(d^{\prime}))}\pure(d^{\prime})_{i-1,j}.
$$
of Lemma~\ref{stepUp} is then
$$
D = \begin{pmatrix} 1&-&-&-\\-&3&2&- \\ -&-&-&- \end{pmatrix} +
\frac{3}{4}\begin{pmatrix} -&1&-&- \\ -&-&2&- \\-&-&-&1 \end{pmatrix}
= \begin{pmatrix} 1&3/4&-&- \\ -&3&7/2&- \\-&-&-&3/4 \end{pmatrix}
$$
and the non-negative linear expansion of $D$ into pure diagrams whose
existence is guaranteed by Lemma~\ref{stepUp} is
$$
D = \frac{2}{5}\pure(0, 1, 3, 5) + \frac{3}{5}\pure(0, 2, 3, 5).
$$
\end{example}

\subsection{Complete intersections}
Let $I$ be an ideal in $R$ such that $R/I$ is a complete intersection, that is, 
$I = (f_1, f_2, \dots, f_{\pdim})$ where $f_1,f_2,\dots,f_{\pdim}$ is a regular sequence
in $R$. 

\begin{thm}\label{CIThm}
The Betti diagram of a complete intersection, $R/I$,
is a non-negative linear combination of pure diagrams.
\end{thm}

\begin{proof}
Denote by $I^{\prime}$ the ideal $(f_1,f_2,\dots,f_{\pdim-1})$.
Then it follows from the Koszul resolutions of $R/I$ and $R/I{^\prime}$
that
\begin{equation} \label{ci_betti}
\beta_{i,j}(R/I) = \beta_{i,j}(R/I{^\prime}) + \beta_{i-1,j-d}(R/I{^\prime})
\end{equation}
where $d$ is the degree of $f_{\pdim}$.
By induction, on the codimension, we may assume that
$$
\beta(R/I{^\prime}) = \sum_{k=1}^d a_k\pure(d_k)
$$
for some pure diagrams $\pure(d_1), \pure(d_2), \dots,\pure(d_{s})$
and some non-negative rational numbers $a_1,a_2,\dots,a_{s}$.
By (\ref{ci_betti}) we get 
$$
\beta(R/I) = \sum_{i=k}^s a_k\rho_k
$$
where $\rho_k$ is defined by
$$
(\rho_k)_{i,j} = \pure(d_k)_{i,j} + \pure(d_k)_{i-1,j-d}.
$$
Hence, to prove that $\beta(R/I)$ is a non-negative linear combination of pure diagrams
it is enough to prove that $\rho_k$ is, for $k = 1,2,\dots,s$.
Assume that $d_k = (d_0,d_1,\dots,d_p)$. Then
$\pure(d_k)_{i-1,j-d} = \pure(d+d_0,d+d_1,\dots,d+d_p)_{i-1,j}$
and we get
$$
(\rho_k)_{i,j} = \pure(d_0,d_1,\dots,d_p)_{i,j}
+ \pure(d+d_0,d+d_1,\dots,d+d_p)_{i-1,j}.
$$
Since the two pure diagrams in the sum giving $\rho_k$ both have
multiplicty $d_1d_2 \cdots d_p$ and clearly $d_i < d + d_i$ for
each $i=0,1,\dots,p$, we get, by
Lemma \ref{stepUp}, that $\rho_k$ is a non-negative linear combination of pure diagrams.
\end{proof}

\subsection{Gorenstein algebras of codimension three}

\begin{thm}\label{GorensteinCodim3Thm}
  Let $A = R/I$ be a Gorenstein algebra of codimension three for some
  ideal $I$ in $R$. Then the Betti diagram of $A$ is a non-negative
  linear combination of pure diagrams, all from the same chain.
\end{thm}

\begin{proof}
  Let $I$ be a height three Gorenstein ideal in $R$ and let
$$
0 \to R(-f) \to \bigoplus_{i=1}^{2k+1}R(-f+a_i) \to
\bigoplus_{i=1}^{2k+1}R(-a_i) \to R \to R/I \to 0
$$
be a minimal free resolution, where $a_1 \leq \dots \leq a_{2k+1}$.
From Geramita and Migliore \cite{Geramita-Migliore} we know that there
is a height two ideal, $I_C$, with minimal free resolution
$$
0 \to \bigoplus_{i=k+2}^{2k+1}R(-f+a_i) \to
\bigoplus_{i=1}^{k+1}R(-a_i)\to R \to R/I_C \to 0.
$$
Let $F_0 = R$, $F_1 = \bigoplus_{i=1}^{k+1}R(-a_i)$ and $F_2 =
\bigoplus_{i=k+2}^{2k+1}R(-f+a_i)$.  Then, numerically, we can write
the resolution of $R/I$ as
$$
\begin{CD} 
  0 @>>> F_0(-f) @>>> F_1(-f) @>>> F_2(-f)  \\
  & & \oplus   & &\oplus   & & \oplus \\
  && 0 @>>> F_2 @>>> F_1 @>>>  F_0 \\
\end{CD}
$$

From Theorem~3.3 we know that the Betti diagram of the resolution $0
\to F_2 \to F_1 \to F_0$ is a non-negative linear combination of pure
diagrams.  These pure diagrams then have shifts $(0,a_i,f-a_j)$, for
some $1 \leq i \leq k+1$ and and $k+2 \leq j \leq 2k+1$, such that
$a_i \leq f-a_j$.  Hence the Betti diagram of $0 \to F_2 \to F_1 \to
F_0$ is
$$
\sum_{i,j} b_{i,j} \pure(0,a_i,f-a_j).
$$
for some non-negative rational numbers $b_{i,j}$.

Let $\pi = \pure(0,a_i,f-a_j)$ be such a diagram and define a diagram
$D(i,j)$ by
$$
D(i,j)_{k,l} = \pure(0,a_i,f-a_j)_{k,l}
+\frac{a_i}{a_i+a_j-f}\pure(a_j,f-a_i,f)_{k-1,l}.
$$
The Betti diagram of $R/I$ is then
$$
\sum_{i,j} b_{i,j}D(i,j),
$$
so all we have to do is to show that $D(i,j)$ is a non-negative linear
combination of pure diagrams.  By Lemma \ref{stepUp}, $D(i,j)$ is a
non-negative linear combination of pure diagrams if the two diagrams
in the sum giving $D(i,j)$ have the same multiplicity and $0 < a_j$,
$a_i < f - a_i$ and $f-a_j < f$. Both diagrams have multiplicity
$a_i(f-a_j)$ so the only thing left to check here is that $a_i < f -
a_i$ for $i=1,\dots,k+1$.  Let $r_i = f - 2a_i$ and note that $r_1
\geq \dots \geq r_{2k+1}$.  We know from \cite[Proposition
3.1]{Diesel} that
$$
r_i+r_{2k+3-i} > 0, \text{ for } i=2,\dots,k+1.
$$
Now we have
$$
r_{k+1} - r_{k+2} \geq 0
$$
and
$$
r_{k+1}+r_{k+2} > 0
$$
and we get
$$
r_1 \geq \dots \geq r_{k+1} > 0
$$
which then implies $a_i < f - a_i$ for $i=1,\dots,k+1$.

It remains to prove that the pure diagrams in the non-negative linear
combination can be taken all from the same chain. The minimal and
maximal shifts of $A$ are $\down = (0,a_1,f-a_{2k+1},f)$ and $\up =
(0,a_{2k+1},f-a_1,f)$.  Since $\up_3=\down_3 =f$, Lemma~\ref{IsoLemma}
gives an isomorphism $\phi_3:\vecsp \to V_{\tau_3(\down),\tau_3(\up)}$
where $\tau_3(\down) = (0,a_1,f-a_{2k+1})$ and $\tau_3(\up) =
(0,a_{2k+1},f-a_1)$.  Since $\beta(A)$ is a non-negative linear
combination of pure diagrams, the same is true for $\phi_3(\beta(A))$
since $\phi_3$ preserves non-negative linear combinations of pure
diagrams. From Theorem~\ref{CodimOneTwo} we get, since
$\phi_3(\beta(A))$ has codimension two, that the pure diagrams in the
non-negative linear combination giving $\phi_3(\beta(A))$ can be
chosen from the same chain.  Since $\phi_3$ preserves chains of pure
diagrams we get that $\beta(A)$ is a non-negative linear combination
of pure diagram all from the same chain.
\end{proof}

\section*{Acknowledgements}
We would like to thank Fabrizio Zanello for introducing us
to the Multiplicity conjecture.

\bibliographystyle{abbrv} \bibliography{/afs/kth.se/home/j/o/jonasso/artlev/latex/algbib}

\end{document}